\documentclass[12pt]{article}
\usepackage[utf8]{inputenc}
\usepackage[english]{babel}
\usepackage{citehack,amsmath,amsthm,amssymb}
\usepackage{psfrag}
\usepackage[active]{srcltx}
\usepackage{graphicx}
\usepackage{epsfig}
\usepackage{enumerate}
\usepackage[usenames]{color}
\usepackage[a4paper,left=28mm,right=12mm,top=2cm,bottom=2cm]{geometry}

\usepackage{xspace,cite,indentfirst}
\usepackage[mathscr]{eucal}

\begin{document}

\begin{center}
{\bf Weak chaos in discrete systems}
\end{center}

~~~~~~~~~~~~~~~~~~~~~~~~~~~~~~~~~~~~~~~~~~~~~Pokutnyi Oleksandr


{\bf Abstract.} For a general discrete dynamics on a Banach and Hilbert spaces we give necessary and
sufficient conditions of the existence of bounded solutions under assumption that the homogeneous difference equation
admits an discrete dichotomy on the semi-axes. We consider the so called resonance (critical) case when the uniqueness of solution is disturbed. We show that admissibility  can be reformulated in the terms of generalized or pseudo-invertibility. As an application we consider the case when the corresponding dynamical system is e-trichotomy.




{\bf Introduction.} The diffeomorphism of $\mathbb{R}^{n}$ into itself with a homoclinic point is chaotic (see \cite{Smale}, \cite{Smale1}).
Using the Melnikov method \cite{Meln} or the theory of Noetherian operators \cite{Lin}, planar periodic systems with the Poincare map were constructed, which have transversal homoclinic points and, therefore, chaotic. Shilnikov examined an autonomous system in space with an orbit homoclinic up to the saddle-node and showed, under certain conditions, the presence of chaos. In the article \cite{Shiln} Shilnikov again considered autonomous systems in space with orbits homoclinic to hyperbolic periodic orbits and showed, under transversality conditions, the existence of chaos. For finite-dimensional systems and periodic mappings, Palmer rigorously proved Smale's horseshoe theorem in 1988 using the shadowing lemma (see \cite{Palm}). For infinite-dimensional mappings and periodic systems, the proof was completed by Steinlein and Walther (see \cite{Stein}). For finite-dimensional autonomous systems, such a proof was performed by Palmer in 1996. For infinite-dimensional systems in \cite{Yang}. It should also be noted recent monographs on chaos  \cite{Zar} (see also \cite{Kalosh}). Proposed in the article method gives possibility to observe so called "weak"$~$ chaos in discrete systems. A separate article will be devoted to this question.
The concept of dichotomy plays a key role in the qualitative
theory of dynamical systems. It is worth to mention here the papers \cite{Bar_Valls,  BarVal,  BarVal1, Bar2, Ben, BoiSam1, Chuesh, Hen, Mee, Meg, Sasu_L_B, Sasu1, Slyus,  Karabash} where dichotomy for difference equations was studied. Exponential-dichotomous systems represent themselves a class of trajectories which can either grow exponentially or decrease with the exponential rate at infinity. 
It generalizes the well-known notion of hyperbolicity of an operator to the non-autonomous case (see \cite{Chuesh}). 
A lot of scientific papers study difference equations that allow an dichotomy on the whole integer axis $\mathbb{Z}$ or on the semi-axes $\mathbb{Z}_{+}$. It should be noted that the condition of the  dichotomy on $\mathbb{Z}$ is usually equivalent to the regularity of the respective system (it
guarantees the existence of solution condition which is unique and bounded on the entire axis
$\mathbb{Z}$). In contrast, the dichotomy condition on the semi-axes can guarantee only normal or strong generalised normal
solvability. Such class of problems belongs to the irregular (resonant) class. Solutions to this
kind of problems do not exist with all heterogeneity, and their number can be infinite. In the
finite-dimensional case sufficiently full study of systems, that allow exponential dichotomy
on the semi-axes, has been done in the book \cite{BoiSam1} where the Noetherian of the
respective problem was shown. Note that the theory of generalized inverse and Moore-Penrose pseudo-inverse operators was used while studying the underlying problem in irregular case (see \cite{BoiSam1, Moor, Penr}). In this paper, we study this class of systems in infinite-dimensional
spaces. It is shown that in contrast to the finite-dimensional case, for the generated operator of
the corresponding equation, there are possibly much more options. Namely, the generated
operator can be strong generalized normally solvable,  normally solvable, $n$-normal, $d$-normal, Noetherian or Fredholm.

{\bf Statement of the problem.} 
Consider the following equation
\begin{equation}  x_{n+1} = A_{n}x_{n} + h_{n}, n \in
	\mathbb{Z}, \label{6.2.1} \end{equation}  
with boundary condition
\begin{equation} lx_{\cdot} = \alpha, \label{6.2.10}
\end{equation}
where  $A_{n} : B \rightarrow B$ - is a set of bounded operators, from the Banach space $B$ into itself. Assume that
\begin{align} A = (A_{n})_{n \in \mathbb{Z}} \in l_{\infty}(\mathbb{Z}, \mathcal{L}(B)), h =
(h_{n})_{n \in \mathbb{Z}} \in l_{\infty}(\mathbb{Z}, B). \nonumber
\end{align}
It means that

\begin{align}|||A||| = \sup_{n \in \mathbb{Z}} ||A_{n}|| < + \infty, |||h||| =
\sup_{n \in \mathbb{Z}}||h_{n}|| < + \infty, \nonumber 
\end{align}
\begin{align}
l : l_{\infty}(\mathbb{Z}, B) \rightarrow B_{1} \nonumber 
\end{align}
is a linear and bounded operator which translates bounded solutions of (\ref{6.2.1}) into the Banach space $B_{1}$, $\alpha$ is an element of Banach space $B_{1}$.
First we formulate the conditions for the existence of bounded solutions of the equation (\ref{6.2.1}).
The corresponding homogeneous difference equation has the following form:
\begin{equation} \label{6.2.2} x_{n + 1} = A_{n}x_{n}. \end{equation} It should be noted that an arbitrary solution of a homogeneous equation can be represented as: $x_{m} =
\Phi(m, n)x_{n}, m \geq n,$ where
\begin{align}\Phi(m, n) = \left\{ \begin{array}{ccccc} A_{m
		-1}A_{m-2}...A_{n}, ~\mbox{if}~ m >n \\
	I, ~\mbox{if}~ m = n. \end{array} \right. \nonumber
\end{align}
It is clear, that $\Phi(m, 0) = A_{m -1}A_{m -2}...A_{0}$. Also, we denote  $U(m): = \Phi(m, 0)$ and $U(0) = I$.

Traditionally \cite{Chuesh}, the mapping  $\Phi(m,n)$ is called the {\bf evolutionary operator} of the problem (\ref{6.2.2}).
Suppose that the equation (\ref{6.2.2}) is exponentially dichotomous on the semiaxes $\mathbb{Z}_{+}$ and $\mathbb{Z}_{-}$ with projectors $P$ and $Q$ in the space $B$ respectively, which means that there are projectors $P (P^{2} = P)$ and $Q (Q^{2} = Q)$, constants $k_{1, 2} \geq 1,$ $0 < \lambda_{1, 2} < 1$ such that
\begin{align}
||U(n)PU^{-1}(m)|| \leq k_{1}\lambda_{1}^{n - m}, n \geq m, \nonumber \\ 
||U(n)(I - P)U^{-1}(m)|| \leq k_{1}\lambda_{1}^{m - n}, m \geq n \nonumber
\end{align}
for arbitrary $m, n \in \mathbb{Z}_{+}$ (dichotomy on $\mathbb{Z}_{+}$);
\begin{align}
||U(n)QU^{-1}(m)|| \leq k_{2}\lambda_{2}^{n - m}, n \geq m \nonumber \\
||U(n)(I - Q)U^{-1}(m)|| \leq k_{2}\lambda_{2}^{m - n}, m \geq n \nonumber
\end{align}
for arbitrary $m, n \in \mathbb{Z}_{-}$ (dichotomy on $\mathbb{Z}_{-}$).

We give the conditions on the existence of solutions of the equation (\ref{6.2.1}) that are bounded on the entire axis under the assumption, that the homogeneous equation (\ref{6.2.2}) admits the exponential dichotomy on the semi-axes $\mathbb{Z}_{+}$, $\mathbb{Z}_{-}$ with the projectors $P$ and $Q$ introduced above. We use this definition for simplicity. The main result can be proved in more general case (see \cite{Hen}) in the same way.

Recall some facts from the theory of generalized inverse and Moore-Penrose pseudo-inverse operators (see \cite{BoiSam1, Goh, kortur, kreinSG}) which we will use a little bit later.

{\bf Definition 1.} (see \cite{BoiSam1}). {\it A linear bounded operator $D$ from one Banach space $B_1$ into another
	Banach space $B_2$ is called {\bf normally-solvable}, if its image $ImD = R(D)$ is closed ($\overline{R(D)} = R(D)$).
}

{\bf Definition 2.} (see \cite{BoiSam1, kortur}). {\it  Operator $D \in \mathcal{L}(B_{1}, B_{2})$ is called {\bf generalized invertible}, if there exists operator denoted by $D^{-} \in \mathcal{L}(B_{1}, B_{2})$, such that
	\begin{align} D^{-}DD^{-} = D^{-}, DD^{-}D = D. \nonumber
	\end{align}
}

\bigskip
Note, that if the operator $D$ has inverse $D^{-1}$, left inverse $D_{l}^{-1}$ or right inverse $D_{r}^{-1}$ then it
coincides with $D^{-}$: $D^{-1} = D^{-}, D_{l}^{-1} = D^{-}, D_{r}^{-1} = D^{-}$ (see \cite{BoiSam1,  Goh,  kortur, kreinSG}).

The following criterion of generalized invertibility is well-known (see \cite{BoiSam1, kortur}).

{\bf Theorem 1.} \label{th1:1} {\it  Operator $D$ has a generalized inverse if and only if the following conditions hold:
	
	i) $D$ is a normally-solvable operator;
	
	ii) The kernel $kerD = N(D)$ of the operator $D$ is a complemented subspace of $B_{1}$ (there exists a subspace $X$ of the Banach space $B_{1}$ such that $N(D) \oplus X = B_{1}$);
	 
	iii) The image $ImD = R(D)$ of the operator $D$ is a complemented subspace of $B_{2}$ (there exists a subspace $Y$ of the Banach space $B_{2}$ such that $R(D) \oplus Y = B_{2}$).
}

If we consider the operator $D \in \mathcal{L}(H_{1}, H_{2})$ from one Hilbert space $H_{1}$ into another Hilbert space $H_{2}$, then from the set of all generalized inverse operators we can choose the operator $D^{+}$ which has the following properties:
\begin{align}
	D^{+}DD^{+} = D^{+},~~~~~~ DD^{+}D = D, \nonumber
\end{align}
\begin{align}
	(DD^{+})^{*} = DD^{+}, ~~~~~~(D^{+}D)^{*} = D^{+}D. \nonumber
\end{align}
The operator $D^{+} \in \mathcal{L}(H_{2}, H_{1})$ is called the {\bf Moore-Penrose pseudo-inverse} to the operator $D$ (see \cite{Moor,Penr}). It follows from the Theorem 1 that the operator $D$ in a Hilbert space has Moore-Penrose pseudo-inverse if and only if it is normally-solvable (which means that the set of its values is closed) (see \cite{BoiSam1, Goh}).

{\bf Main Results. Linear case.} 

{\bf Banach space case.}
We state and prove the main result of the existence of bounded solutions of the nonhomogeneous equation  (\ref{6.2.1}) defined in Banach space.

{\bf Theorem 2.} \label{th6:20}  {\it Suppose that a homogeneous equation is exponential dichotomous on the semi-axes $\mathbb{Z}_{+}$ and $\mathbb{Z}_{-}$ with projectors $P$ and $Q$ respectively, and the operator
	$$D = P - (I - Q): B \rightarrow B$$
	is generalised invertible. The solutions of the equation (\ref{6.2.1}) bounded on the entire axis $\mathbb{Z}$ exist if and only if the following condition is satisfied
	\begin{equation}
		\label{6.2.3} \sum_{k = -\infty}^{+\infty}H(k+1)h_{k} = 0.
	\end{equation}
	If the condition (\ref{6.2.3}) holds, the set of bounded solutions has the following view:
	\begin{equation} \label{6.2.40}
		x_{n}(c) = U(n)P\mathcal{P}_{N(D)}c + (G[h])(n),~~~~~c \in B
	\end{equation}
	where
	\begin{align}
	G[h](n)  = \left\{ \begin{array}{ccccc}
		\sum_{k = 0}^{n - 1}U(n)PU^{-1}(k + 1)h_{k} - \sum_{k = n}^{+
			\infty}U(n)(I - P)U^{-1}(k + 1)h_{k} + \nonumber \\
		+ U(n)PD^{-}[\sum_{k = 0}^{+\infty}(I - P)U^{-1}(k + 1)h_{k}+ \nonumber \\
		+ \sum_{k = -\infty}^{-1}QU^{-1}(k + 1)h_{k}],
		\hspace{0.2cm} n \geq 0 \nonumber \\
		\sum_{k = -\infty}^{n - 1}U(n)QU^{-1}(k + 1)h_{k} - \sum_{k =
			n}^{-1}U(n)(I - Q)U^{-1}(k + 1)h_{k} + \nonumber \\
		+ U(n)(I - Q)D^{-}[\sum_{k = 0}^{+ \infty}(I - P)U^{-1}(k + 1)h_{k} + \nonumber \\
		+ \sum_{k = -\infty}^{-1}QU^{-1}(k + 1)h_{k}], \hspace{0.2cm} n
		\leq 0, \nonumber
     \end{array}
	\right. 
    \end{align}
	is generalised inverse Green's operator on $\mathbb{Z}$ with the following properties:
	\begin{align}
	(G[h])(0 + 0) - (G[h])(0 - 0) = - \sum_{k = -\infty}^{+ \infty} H(k + 1)h_{k} = 0, \nonumber
	\end{align}
	\begin{align}
	(LG[h])(n) = h_{n}, n \in \mathbb{Z}, \nonumber
	\end{align}
	where
	\begin{align}
	(Lx)(n):= x_{n + 1} - A_{n}x_{n} : l_{\infty}(\mathbb{Z}, B)
	\rightarrow l_{\infty}(\mathbb{Z}, B), \nonumber
	\end{align}
	$H(k + 1) = \mathcal{P}_{B_{D}}QU^{-1}(k + 1) = \mathcal{P}_{B_{D}}(I -P)U^{-1}(k + 1), D^{-}$  is generalized inverse to the operator $D$, projectors $\mathcal{P}_{N(D)} = I - D^{-}D$ and  $\mathcal{P}_{B_{D}} = I - DD^{-}$ (see \cite{Goh, BoiSam1}), which project space $B$
	on the kernel $N(D)$ of the operator $D$ and the subspace $B_{D}: = B \ominus R(D)$  respectively ($B = B_{D} \oplus R(D) = (B \ominus R(D)) \oplus R(D)$).
}

{\bf Proof.} The proof follows the same scheme as in the book (\cite[p. 288-297]{BoiSam1}). 

{\bf Remark 1.}{\it We have the following estimates for the norm of the solution
	\begin{align}
	||x_{n}|| \leq k_1\lambda_{1}^{n}||\mathcal{P}_{N(D)}c|| + k_1 \lambda_1^{n}||D^{-}|| \left(\frac{k_1\lambda_1}{1 - \lambda_1} + \frac{k_2\lambda_2}{1 - \lambda_2}\right)|||h||| +  \nonumber \\
	+ k_1\frac{(1 + \lambda_1 - \lambda_1^{n})}{1 - \lambda_1}|||h|||, ~~ n  \geq 0, \nonumber
	\end{align}
	\begin{align}
	||x_{n}|| \leq k_2 \lambda_2^{-n}||\mathcal{P}_{N(D)}c|| + k_2 \lambda_{2}^{-n}||D^{-}|| \left(\frac{k_1\lambda_1}{1 - \lambda_1} + \frac{k_2\lambda_2}{1 - \lambda_2}\right)|||h||| + \nonumber \\
	+ k_2 \frac{1 + \lambda_2 - \lambda_2^{-n + 1}}{1 - \lambda_2}|||h|||,~~ n \leq 0. \nonumber
	\end{align}
	From these inequalities follows the estimate
	\begin{align}
	|||x||| \leq K ||\mathcal{P}_{N(D)}c|| +  K||D^{-}|| \left(\frac{k_1\lambda_1}{1 - \lambda_1} + \frac{k_2\lambda_2}{1 - \lambda_2}\right)|||h||| + K\frac{1 + \Lambda_1}{1 - \Lambda_2}|||h|||, \nonumber
	\end{align}
	where $K = \max\{ k_1, k_2\}$, $\Lambda_1 = \max\{ \lambda_1, \lambda_2\}$, $\Lambda_2 = \min\{ \lambda_1, \lambda_2 \}$.
}

{\bf Remark 2.}{\it It should also be noted that if bounded solutions are united together at zero as follows
	\begin{align}
		x_{0+} = x_{0-} + c, \nonumber
	\end{align}
	where $c$ is an element of Banach space, then we obtain a bounded solution for problem with a jump.
}

{\bf Remark 3.}{\it It should be noted that the condition of solvability (\ref{6.2.3})        
    \begin{align}
		\mathcal{P}_{B_{D}}g= 0, \nonumber
    \end{align}
	of theorem \ref{th6:20} is equivalent to
	\begin{align}
		\mathcal{P}_{N(D^{*})}g = 0, \nonumber
	\end{align}
	where $\mathcal{P}_{N(D^{*})}$ is projector onto the kernel of the operator $D^{*}$ adjoint to $D$.
	We also have that $ {}^\bot N(D^{*}) = \{ x \in B : \phi(x) = 0, \phi \in N(D^{*}) \} = R(D)$ (see \cite{kreinSG}).
}

 {\bf Remark 4.} {\it If the operator $D$ is topologically Noetherian (see \cite{AbdBur}) then it has generalized invertible. If the operator $D$ is $n$-normal (see \cite{kreinSG}) then in lemma 1 projectors $\mathcal{P}_{N(D)}$ and $P\mathcal{P}_{N(D)}$ are finite-dimensional and equation (\ref{6.2.1}) has a finite family of bounded solutions (its number depends on the dimension of the operator $\mathcal{P}_{B_{D}}Q$). In the case when $D$ is $d$-normal we have a finite-dimensional projectors $\mathcal{P}_{B_{D}}$, $\mathcal{P}_{B_{D}}Q$ and finite number conditions of solvability in (\ref{6.2.3}) (its number depends on the dimension of the operator $\mathcal{P}_{B_{D}}Q$). If the operator $D$ is Noetherian we obtain that the considering equation (\ref{6.2.1}) has a finite nonzero index, finite number conditions of solvability and finite family of bounded solutions ($ind~L = dim~P\mathcal{P}_{N(D)} - dim~\mathcal{P}_{B_{D}}Q$).
}

{\bf Remark 5.} {\it Sometimes (see  \cite{BoiSam1}, \cite{kreinSG}), Fredholm operators are called Noetherian operators. The Noetherian operators are operators for which all conditions of Fredholm with zero index operator are satisfied but $N(L)$ may differ from $m(L)$. This class of operators is named after F. Noether who studied a class of singular integral equations with operators of this sort for the first time as early as in 1920 (see \cite{Noether}).}

{\bf Remark 6.} {\it Suppose that under condition of Theorem \ref{th6:20}
$$[P,Q] = PQ - QP = 0,~~ PQ = Q.$$ This is the so-called case of exponential trichotomy of equation
(\ref{6.2.2}) on  $\mathbb{Z}$. In this case, inhomogeneous equation (\ref{6.2.1}) has at least one solution on $\mathbb{Z}$ for any element $h \in l_{\infty}(\mathbb{Z}, B)$.
}

{\bf Theorem 2.} {\it If the conditions of the Theorem \ref{th6:20} and additional conditions such as following hold $$[P,Q] =PQ - QP = 0, PQ = Q = P$$ then the nonhomogeneous equation (\ref{6.2.1}) has unique bounded on $\mathbb{Z}$ solution for arbitrary
	$ h \in l_{\infty}(\mathbb{Z}, B)$. }

{\bf Remark 7.}{\it In this case, the considered system is exponentially-dichotomous on the entire axis $\mathbb{Z}$.
	In the finite dimensional case, this result is well known (see \cite{Chuesh}). Theorem \ref{th6:20} under less restrictive assumptions, allows to find a set of bounded solutions. }

Let us find the solvability condition of the boundary-value problem (\ref{6.2.1}), (\ref{6.2.10}). Suppose that condition (\ref{6.2.3}) is satisfied. Substitute expression (\ref{6.2.40}) in the boundary condition (\ref{6.2.10}). Then we obtain the following operator equation 
\begin{equation} \label{equation:1}
Vc = \alpha - l(G[h])(\cdot),
\end{equation}
where $V = lU(\cdot)P\mathcal{P}_{N(D)} : B \rightarrow B_{1}$. If the operator $V$ is generalized invertible then equation (\ref{equation:1}) is solvable if and only if the following condition is satisfied
\begin{equation} \label{equation:2}
\mathcal{P}_{V_{B_{1}}}\left( \alpha - l(G[h])(\cdot)\right) = 0,
\end{equation}
where $\mathcal{P}_{V_{B_{1}}} = I - VV^{-} (B_{1} = R(V) \oplus V_{B_{1}})$. Under condition (\ref{equation:2}) the set of solutions of the equation (\ref{equation:1}) has the following form 
\begin{align}
c = V^{-}\left(\alpha - l(G[h])(\cdot) \right) + \mathcal{P}_{N(V)}\overline{c}, \overline{c} \in B. \nonumber
\end{align}
Thus, we obtain the following theorem. 

{\bf Remark 8.} {\it It should be noted that the operator $l$ in (\ref{6.2.10}) can be for example two-point or multi-point 
\begin{align}
lx_{\cdot} = A_{1}x_{m} - A_{2}x_{0}, ~~~lx_{\cdot} = \sum_{i = 1}^{n}A_{i}x_{m_{i}}, \nonumber
\end{align}
where $m, m_{i} \in \mathbb{Z}$, $A_{i} \in \mathcal{L}(B, B_{1}), i = \overline{1, n}$ are linear and bounded operators. Another example is conditions at the infinity 
\begin{align}
lx_{\cdot} = A_{1} \lim_{n \rightarrow -\infty}x_{n} + A_{2} \lim_{n \rightarrow +\infty}x_{n}. \nonumber
\end{align}
}

{\bf The case of Hilbert spaces.} 

Let us show, that in the case of Hilbert spaces it is possible to obtain a more general result. Using the technique, given in \cite{BoiPok100} it is possible to introduce the concept of a strong pseudo-inverse operator to an operator  $D$ and weaken the assumption of the Theorem \ref{th6:20} on the generalized invertibility of the operator $D$.
In this case, the original equation (\ref{6.2.1}) can always be made solvable in one of the meanings, defined below. Let us explain this approach in more details.

Let us describe the construction of a strong pseudo-inverse according to Moore-Penrose
operator, which is used to represent the solutions
of the operator equation of the following operator equation
\begin{equation} \label{news:1}
	D\xi = g, ~~~D: \mathcal{H} \rightarrow \mathcal{H},
\end{equation}
$\mathcal{H}$ - is a Hilbert space, $g \in \mathcal{H}$.
We distinguish three types of solutions.

1) Classical solutions.

If the operator $D$ is normally solvable $(R(D)=\overline{R(D)})$ and, therefore, pseudo-invertible, then the element $g$ belongs to the set of values  $(g \in R(D))$ of the operator $D$ if and only if $P_{B_{D}}g = 0$ or $P_{N(D^{*})}g = 0$ (see \cite{BoiSam1}); $P_{B_{D}}$, $P_{N(D^{*})}$ --
are orthoprojectors onto the subspace $B_{D}$ and cokernel of the operator $D$. In this
case, there is Moore-Penrose pseudo-inverse operator $D^{+}$ and the set of solutions of the
equation (\ref{news:1}) has the form
\begin{align}
	\xi = D^{+}g + P_{N(D)}c, ~\mbox{for any}~ c \in \mathcal{H}, \nonumber
\end{align}
where $P_{N(D)}$ - is orthoprojector onto the kernel of the operator $D$.

2) Strong generalized solutions.

Consider a case, when the set of values of the operator  $D$ is not closed, i.e. $R(D) \neq \overline{R(D)}$. We
show that in this case $D$ can be extended to the operator $\overline{D}$ in the way that the operator $\overline{D}$ is normally-solvable.

Since the operator $D$  is linear and bounded, the decompositions of the space $\mathcal{H}$ into direct
sums take place
\begin{align}
	\mathcal{H} = N(D) \oplus X, \mathcal{H} = \overline{R(D)}
	\oplus Y. \nonumber 
\end{align}
Here $X = N(D)^{\bot}$,  $Y = \overline{R(D)}^{\bot}$. It can be stated, that there are operators of orthogonal projection $P_{N(D)},~ P_{X}$ onto the corresponding spaces
$P_{\overline{R(D)}}, ~ P_{Y}$. Denote by
$\mathcal{H}_{1}$ the quotient space of the space $\mathcal{H}$ with respect to the
kernel $N(D)$ ($\mathcal{H}_{1} = \mathcal{H}/N(D)$). Then, as it is known (see \cite{Atiyah,Rob}) there exist a continuous bijection $p : X \rightarrow \mathcal{H}_{1}$ and projection $j : \mathcal{H}
\rightarrow \mathcal{H}_{1}$. The triple  $(\mathcal{H}, \mathcal{H}_{1}, j)$ is a locally trivial bundle with a typical layer $
P_{N(D)}\mathcal{H}$. Now we
define the operator
\begin{align}
	\mathcal{D} = P_{\overline{R(D)}}Dj^{-1}p : X
	\rightarrow R(D) \subset \overline{R(D)}. \nonumber
\end{align}
It is easy to verify that the operator which is defined in this way is linear, injective and
continuous. Now using the completion process (see \cite{Lyash}) by norm
$||x||_{\overline{X}} = ||\mathcal{D}x||_{F}$, where $F =
\overline{R(D)}$, we obtain a
new space $\overline{X}$ and an extended operator $\overline{\mathcal{D}}$. Then
\begin{align}
	\overline{\mathcal{D}} : \overline{X} \rightarrow \overline{R(D)},
	~~X \subset \overline{X} \nonumber
\end{align}
and the operator constructed in this way is a homeomorphism between the spaces $\overline{X}$ and $\overline{R(D)}$.
Consider the extended operator $\overline{D} =
\overline{\mathcal{D}}P_{\overline{X}} : \overline{\mathcal{H}}
\rightarrow \mathcal{H}$,
\begin{align}
	\overline{\mathcal{H}} = N(D) \oplus \overline{X},  \mathcal{H} = R(\overline{D})
	\oplus Y. \nonumber
\end{align}
It is clear that $\overline{D}x = Dx, x \in \mathcal{H}$ and operator
$\overline{D}$ is normally-solvable (in this case $R(\overline{D}) =
\overline{R(\overline{D})}$), therefore
is pseudo-invertible with pseudo-inverse $\overline{D}^{+}$ (see \cite{BoiSam1}).

{\bf Definition 2.} {\it The operator $\overline{D}^{+}$ is called {\bf strong pseudo-inverse} to the operator $D$.}

Then the set of strong generalized solutions of the equation  (\ref{news:1}) has the form:
\begin{align}
	\xi = \overline{D}^{+}g + P_{N(\overline{D})}c, ~~\mbox{for any}~ c \in \overline{\mathcal{H}}. \nonumber
\end{align}

Note also, that if in this case $g \in R(D)$, then strong generalized solution becomes the usual classic
solution.

3) Strong generalized pseudo-solutions.

Let us consider the case, when $g \notin \overline{R(D)}$, which for the element $g$ is equivalent to the condition $P_{N(D^{*})}g \neq 0$. In this case, there are elements of $\overline{\mathcal{H}}$,
which minimize the norm $||\overline{D}\xi - g ||_{\mathcal{H}}$ for
$\xi \in \overline{\mathcal{H}}$,
\begin{align}
	\xi = \overline{D}^{+}g + P_{N(\overline{D})}c,~~c \in
	\overline{\mathcal{H}}. \nonumber
\end{align}
These elements will be called {\bf generalized pseudo-solutions} of equation (\ref{news:1}) by analogy with
pseudo-solutions (see \cite{BoiSam1}). Note that if $\overline{\mathcal{H}} = \mathcal{H}$ (the operator $D$ has a closed set of values), then
strong generalized pseudo-solutions are actually ordinary pseudo-solutions.

Combining now the above with the results obtained earlier we can formulate the following lemma.

{\bf Theorem 3.} \label{th6:2}  {\it Suppose that a homogeneous equation is exponentially-dichotomous on the semi-axes $\mathbb{Z}_{+}$ and $\mathbb{Z}_{-}$ with projectors $P$ and $Q$ respectively.
	
	1 a). Strong generalized solutions of the equation (\ref{6.2.1}), bounded on the entire axis $\mathbb{Z}$ exist if and only if the following condition is satisfied
	\begin{equation}
		\label{6.20.100} \sum_{k = -\infty}^{+\infty}\overline{H}(k+1)h_{k} = 0.
	\end{equation}

	1 b). Under the condition (\ref{6.20.100}), the set of bounded strong generalized solutions has the following
	form:
	\begin{equation} \label{6.2.4}
		x_{n}^{0}(c) = U(n)PP_{N(D)}c + \overline{(G[h])}(n),~~~~~c \in \overline{\mathcal{H}}
	\end{equation}
	where $\overline{G[h]}(n)$ - is the extension of the generalized Green operator on  $\overline{\mathcal{H}}$,
	\begin{align}
		\overline{G[h]}(n)  = \left\{
		\begin{array}{rcl}
			\sum_{k = 0}^{n - 1}U(n)PU^{-1}(k + 1)h_{k} - \sum_{k = n}^{+
				\infty}U(n)(E - P)U^{-1}(k + 1)h_{k} + \\
			+ U(n)P\overline{D}^{+}[\sum_{k = 0}^{+\infty}(E - P)U^{-1}(k + 1)h_{k}+ \\
			+ \sum_{k = -\infty}^{-1}QU^{-1}(k + 1)h_{k}],
			\hspace{0.2cm} n \geq 0\\
			\sum_{k = -\infty}^{n - 1}U(n)QU^{-1}(k + 1)h_{k} - \sum_{k =
				n}^{-1}U(n)(E - Q)U^{-1}(k + 1)h_{k} + \\
			+ U(n)(E - Q)\overline{D}^{+}[\sum_{k = 0}^{+ \infty}(E - P)U^{-1}(k + 1)h_{k} +\\
			+ \sum_{k = -\infty}^{-1}QU^{-1}(k + 1)h_{k}], \hspace{0.2cm} n
			\leq 0,
		\end{array}
		\right. \nonumber
	\end{align}
	
	2 a). Strong generalized pseudo-solutions of equation (\ref{6.2.1}) bounded on the entire axis $\mathbb{Z}$ exist if and  only if the following condition is satisfied
	\begin{equation}
		\label{6.2.1000} \sum_{k = -\infty}^{+\infty}\overline{H}(k+1)h_{k} \neq 0.
	\end{equation}
	
	2 b). Under the condition (\ref{6.2.1000}), the set of bounded strong generalized pseudo-solutions has the following
	form:
	\begin{align}
		x_{n}^{0}(c) = U(n)PP_{N(D)}c + \overline{(G[h])}(n),~~~~~c \in \overline{\mathcal{H}}, \nonumber
	\end{align}
	where $\overline{H}(n + 1) = P_{\overline{\mathcal{H}}_{\overline{D}}}QU^{-1}(n + 1) = P_{\overline{\mathcal{H}}_{\overline{D}}}(E -
	P)U^{-1}(n + 1), \overline{D}^{+}$  is strong pseudo-inverse according to Moore-Penrose operator to the operator $D$, $P_{N(\overline{D})} = I - \overline{D}^{-}\overline{D}$ and  $P_{\overline{\mathcal{H}}_{\overline{D}}} = I - \overline{D}\overline{D}^{-}$ are projectors that project the space $\mathcal{H}$
	onto the kernel  $N(\overline{D})$ of the operator $\overline{D}$ and the subspace $\overline{\mathcal{H}}_{\overline{D}}: = \overline{\mathcal{H}} \ominus R(\overline{D})$  respectively ($\overline{\mathcal{H}} = \overline{\mathcal{H}}_{\overline{D}} \oplus R(\overline{D}) = (\overline{\mathcal{H}} \ominus R(\overline{D})) \oplus R(\overline{D})$).
}

The following corollaries are also true.

{\bf Theorem 4.} {\it Suppose that $[P, Q] = 0$ and equation (\ref{6.2.1}) is considered in Hilbert space $B = \mathcal{H}$. Then the operator $D$ has a generalized inverse  $D^{-}$,
	which coincides 		with the operator $D$ ($D^{-} = D$). In this situation, we have the following options:
	
	1 a.) Equation (\ref{6.2.1}) has bounded solutions if and only if the following condition is satisfied
	
	\begin{equation}
		\label{6.2.100} \sum_{k = -\infty}^{+\infty}H(k+1)h_{k} = 0;
	\end{equation}
	
	1 b.) When the condition (\ref{6.2.100}) is satisfied, bounded solutions have the form
	
	\begin{equation} \label{6.2.400}
		x_{n}^{0}(c) = U(n)(P - PQ)c + (G[h])(n),
	\end{equation}
	where
	\begin{align}
		G[h](n)  = \left\{
		\begin{array}{rcl}
			\sum_{k = 0}^{n - 1}U(n)PU^{-1}(k + 1)h_{k} - \sum_{k = n}^{+
				\infty}U(n)(E - P)U^{-1}(k + 1)h_{k} + \\
			+ U(n)PD[\sum_{k = 0}^{+\infty}(E - P)U^{-1}(k + 1)h_{k}+ \\
			+ \sum_{k = -\infty}^{-1}QU^{-1}(k + 1)h_{k}],
			\hspace{0.2cm} n \geq 0\\
			\sum_{k = -\infty}^{n - 1}U(n)QU^{-1}(k + 1)h_{k} - \sum_{k =
				n}^{-1}U(n)(E - Q)U^{-1}(k + 1)h_{k} + \\
			+ U(n)(E - Q)D[\sum_{k = 0}^{+ \infty}(E - P)U^{-1}(k + 1)h_{k} +\\
			+ \sum_{k = -\infty}^{-1}QU^{-1}(k + 1)h_{k}], \hspace{0.2cm} n
			\leq 0,
		\end{array}
		\right. \nonumber
	\end{align}
	is generalized Green's operator
	$H(n + 1) = P_{\mathcal{H}_{D}}QU^{-1}(n + 1) = P_{\mathcal{H}_{D}}(E -
	P)U^{-1}(n + 1) = (Q - PQ)U^{-1}(n + 1),$  $P_{N(D)} = I - D^{2} = P_{\mathcal{H}_{D}} $ are projectors.
	
	2 a.) Equation (\ref{6.2.1}) has pseudo-solutions if and only if
	\begin{equation}
		\label{6.2.101} \sum_{k = -\infty}^{+\infty}H(k+1)h_{k} = q \neq 0;
	\end{equation}
	
	2 b.) Under condition (\ref{6.2.101}) the set of bounded pseudo-solutions has this view
	\begin{equation} \label{6.2.401}
		x_{n}^{0}(c) = U(n)(P -PQ)c + (G[h])(n),
	\end{equation}
	where
	\begin{align}
		G[h](n) = \left\{
		\begin{array}{rcl}
			\sum_{k = 0}^{n - 1}U(n)PU^{-1}(k + 1)h_{k} - \sum_{k = n}^{+
				\infty}U(n)(E - P)U^{-1}(k + 1)h_{k} + \\
			+ U(n)PD[\sum_{k = 0}^{+\infty}(E - P)U^{-1}(k + 1)h_{k}+ \\
			+ \sum_{k = -\infty}^{-1}QU^{-1}(k + 1)h_{k}],
			\hspace{0.2cm} n \geq 0\\
			\sum_{k = -\infty}^{n - 1}U(n)QU^{-1}(k + 1)h_{k} - \sum_{k =
				n}^{-1}U(n)(E - Q)U^{-1}(k + 1)h_{k} + \\
			+ U(n)(E - Q)D[\sum_{k = 0}^{+ \infty}(E - P)U^{-1}(k + 1)h_{k} +\\
			+ \sum_{k = -\infty}^{-1}QU^{-1}(k + 1)h_{k}], \hspace{0.2cm} n
			\leq 0,
		\end{array}
		\right. \nonumber
	\end{align}
	is generalized Green's operator.
}

{\bf Proof.} If $[P, Q] = 0$, then it is easy to see that 
\begin{align}
	DDD = D. \nonumber
\end{align}
From this equality, it follows that $D = D^{-}$. 
The second part of the Theorem follows from the fact that the equation (\ref{news:1}) under the condition
\begin{align}
	P_{B_{D}}g  = P_{\mathcal{H}_{D}}g \neq 0, \nonumber
\end{align}
has (see \cite{BoiSam1}) pseudo-solutions (in this case this condition is equivalent to $g \notin R(D)$) and then the set
$\xi = D^{-}g + P_{N(D)}c, c \in \mathcal{H}$ is
the set of those elements which minimize the residual norm $||D\xi - g||$. Substituting the obtained identities
into representations (\ref{6.2.401}) we obtain
\begin{align}
	PP_{N(D)} = P(I - D^{2}) = P(P + Q - 2PQ) = P^{2} + PQ - 2P^{2}Q = P - PQ, \nonumber \\
	P_{B_{D}}Q = P_{\mathcal{H}_{D}}Q = (I - D^{2})Q = (P + Q - 2PQ)Q = \\
	= PQ + Q^{2} - 2PQ^{2} = Q - PQ.
\end{align}

{\bf Theorem 5.}{\it If additional $PQ = 0$ then
	\begin{align}
		PP_{N(D)} = P, P_{B_D}Q = Q \nonumber
	\end{align}
}

{\bf Theorem 6.}	{\it Suppose that $[P, Q] = 0$ and $P = P^{*}$, $Q = Q^{*}$. Then the operator $D$ has a Moore-Penrose pseudo-inverse, which coincides with the operator $D$ ($D^{+} = D$).
}

{\bf Proof.}	This follows from the fact that $DDD = D$ and $D^{*} = P^{*} - I + Q^{*} = D$.

Consider the case when the condition (\ref{6.2.1000}) is satisfied. Using the construction of strong Moore-Penrose operator and theorem 1 we can obtain the following result (in this case we assume that the operator $l: l_{\infty}(\mathbb{Z}, \overline{\mathcal{H}}) \rightarrow \mathcal{H}_{1}$).  

{\bf Theorem 7.} {\it Suppose that condition (\ref{6.2.1000}) is satisfied. Then:

a) Under condition 
\begin{equation} \label{equation:3}
P_{\overline{V}_{\mathcal{H}_{1}}}\left(\alpha - l(G[h])(\cdot) \right) = 0
\end{equation}
boundary-value problem (\ref{6.2.1}), (\ref{6.2.1000}) has a set of strong generalized bounded solutions in the form 
\begin{equation} \label{equation:4}
x_{n}^{0}(\overline{c}) = U(n)PP_{N(D)}P_{N(V)}\overline{c} + \overline{G[h, \alpha]}(n), \overline{c} \in \mathcal{H}
\end{equation}
where 
\begin{align}
\overline{G[h, \alpha]}(n) = (G[h])(n) + \overline{V}^{+} \left( \alpha - l(G[h])(\cdot)\right)
\end{align}
is the extension of the generalized Green's operator on $\widetilde{\mathcal{H}}: \overline{\mathcal{H}} \subset \widetilde{\mathcal{H}} (\widetilde{\mathcal{H}}$ is the completion of $\overline{\mathcal{H}})$ according to the corresponding norm and $P_{\overline{V}_{\mathcal{H}_{1}}, P_{N(V)}}$ are projectors onto the corresponding subspaces. Here $\overline{V} : \widetilde{\mathcal{H}} \rightarrow \mathcal{H}_{1}$, space $\mathcal{H}_{1}$ has the following representation $\mathcal{H}_{1} = R(\overline{V}) \oplus \overline{V}_{\mathcal{H}_{1}}$.

b) Under condition 
\begin{equation} \label{equation:100}
P_{\overline{V}_{\mathcal{H}_{1}}}\left(\alpha -l(G[h])(\cdot) \right) \neq 0
\end{equation}
boundary-value problem (\ref{6.2.1}), (\ref{6.2.1000}) has a set of strong generalized pseudo-solutions in the form
\begin{equation} \label{equation:1000}
x_{n}^{0}(\overline{c})) = U(n)PP_{N(D)}P_{N(V)}\overline{c} + \overline{G[h, \alpha]}(n), \overline{c} \in \mathcal{H}
\end{equation}
}

{\bf Nonlinear case. Bifurcation of solutions.}

Consider the following weakly nonlinear boundary-value problem 
\begin{equation} \label{today:1}
	x_{n + 1}(\varepsilon) = A_{n}x_{n}(\varepsilon) + \varepsilon Z(x_{n}(\varepsilon), n, \varepsilon) + h_{n},
\end{equation}
\begin{equation} \label{today:2}
lx_{\cdot}(\varepsilon) = \alpha
\end{equation}
in the Hilbert space $\mathcal{H}, \mathcal{H}_{1}$ where the nonlinear vector-valued function $Z(x(n, \varepsilon), n, \varepsilon)$ satisfies the following conditions
\begin{align}
Z(\cdot, n, \varepsilon) \in C [||x - x^{0}|| \leq q], Z(x(n, \varepsilon), \cdot, \varepsilon) \in l_{\infty}(\mathbb{Z}, \mathcal{H}), Z(x(n, \varepsilon), n, \cdot) \in C[0, \varepsilon_{0}] \nonumber
\end{align}
in the neighborhood of solution $x_n^{0}(c)$ of the generating linear boundary-value problem (\ref{6.2.1}), (\ref{6.2.1000}) ($q$ is a small enough constant). We are looking for necessary and sufficient conditions for the existence of strong generalized solutions $x_{n}(\varepsilon): \mathbb{Z} \rightarrow \mathcal{H}$ of (\ref{today:1}), (\ref{today:2}) bounded on the entire integer axis
\begin{align}
x_{\cdot}(\varepsilon) \in l_{\infty} (\mathbb{Z}, \mathcal{H}),~~ x_{n}(\cdot) \in C[0, \varepsilon_{0}], \nonumber
\end{align}
which turn into one of the solutions $x_{n}^{0}(c)$ of the generating boundary-value problem (\ref{6.2.1}), (\ref{6.2.1000}) for $\varepsilon = 0$: $x_{n}(0) = x_{n}^{0}(c)$ (in the form (\ref{6.2.4})).

{\bf Theorem 8.} (necessary condition). \label{th:4}
{\it Suppose the equation (\ref{6.2.2}) admits a dichotomy on the semi-axes $\mathbb{Z}_{+}$ and $\mathbb{Z}_{-}$ with projectors $P$ and $Q$ respectively. Let the boundary-value problem (\ref{today:1}), (\ref{today:2}) has a strong generalized solution $x_{n}(\varepsilon)$ bounded on $\mathbb{Z}$, which turns into one of the generating solutions $x_{n}^{0}(\overline{c})$  (\ref{6.2.4}) of the boundary-value problem (\ref{6.2.1}), (\ref{6.2.1000}) with element $\overline{c} = c^{*} \in \overline{\mathcal{H}}$. Then the element $c^{*}$ satisfies the equation
	\begin{equation} \label{today:2}
		F(c^{*}) = \left\{ \begin{array}{ccccc} \sum_{k = -\infty}^{+\infty} \overline{H}(k + 1)Z(U(k)PP_{N(D)}c^{*} + \overline{(G[h])}(k), k, 0) = 0, \\
  P_{\overline{V}_{\mathcal{H}_{1}}}lZ\left(U(\cdot)PP_{N(D)}P_{N(V)}c^{*} + \overline{(G[h, \alpha])}(\cdot), \cdot, 0 \right) = 0.
  \end{array}
  \right.
	\end{equation}
}

{\bf Proof.}	We can consider the nonlinearity $Z$ as an inhomogeneity and apply lemma 2 for the existence of strong generalised bounded solutions of the equation (\ref{6.2.1}) and theorem 2 for the existence of solutions of boundary-value problem (\ref{6.2.1}), (\ref{6.2.1000}). Then we obtain the following conditions of solvability
$$
\sum_{k = -\infty}^{+\infty} \overline{H}(k + 1)Z(x_{k}(\varepsilon), k, \varepsilon) = 0,
$$
and 
$$
P_{\overline{V}_{\mathcal{H}_{1}}}\left(\alpha - l(\overline{G[h + \varepsilon Z(x_{\cdot}, \cdot, \varepsilon), \alpha]})(\cdot) \right) = 0.
$$
Using conditions  (\ref{6.20.100}), (\ref{equation:3}) of strong generalized solvability of linear generating boundary-value problem (\ref{6.2.1}), (\ref{6.2.1000}) and passing to the limit as $\varepsilon \rightarrow 0$ we obtain the required condition of solvability (\ref{today:2}).

{\bf Remark 9.}{\it Equation (\ref{today:2}) will be called the {\bf equation for generating elements} by analogy with the case of periodic problem \cite{BoiSam1}, \cite{Greben}. }

In order to obtain a sufficient condition for the existence of strong generalized bounded solutions, we additionally assume that the nonlinear vector-function $Z$ is strongly differentiable $Z(\cdot, n, \varepsilon) \in C^{1}[||x - x^{0}|| \leq q]$ in a neighborhood of the generating solution $x_{n}^{0}(c^{*})$ (\ref{6.2.4}) ($c^{*}$ is the root of the equation for generating elements (\ref{today:2})).
Making a change of variables $x_{n}(\varepsilon) = y_{n}(\varepsilon) + x_{n}^{0}(c^{*})$ we obtain the following boundary-value problem
\begin{equation} \label{today:3}
	y_{n + 1}(\varepsilon) = A_n y_{n}(\varepsilon) + \varepsilon Z(y_n(\varepsilon) + x_{n}^{0}(c^{*}), n, \varepsilon), ~~y_{n}(0) = 0.
\end{equation}
\begin{equation} \label{today:4}
ly_{\cdot} = 0.
\end{equation}

Due to strong differentiability of $Z$, we obtain the following expansion
\begin{align}
Z(y_{n}(\varepsilon) + x_{n}^{0}(c^{*}), n, \varepsilon) = Z(x_{n}^{0}(c^{*}), n, 0) + A_1(n)y_{n}(\varepsilon) + \mathcal{R}(y_n(\varepsilon), n, \varepsilon)
\end{align}
where 
\begin{align}
A_1(n) = \frac{\partial Z(x, n, 0)}{\partial x}|_{x = x_{n}^{0}(c^{*})} \nonumber 
\end{align}
is Frechet derivative and nonlinearity $\mathcal{R}$ such that $ \mathcal{R}(0, n, 0) = 0$. Under condition of strong generalized solvability
(\ref{6.2.4}) of the equation (\ref{6.2.1}) 
\begin{equation} \label{today:4}
	\sum_{k = -\infty}^{+\infty} \overline{H}(k + 1)\left( Z(x_{k}^{0}(c^{*}), k, 0) + A_1(k)y_{k}(\varepsilon) + \mathcal{R}(y_k(\varepsilon), k, \varepsilon)\right) = 0
\end{equation}
and strong generalised solvability of boundary-value problem (\ref{6.2.1}), (\ref{6.2.1000})
\begin{equation} \label{today:5}
P_{\overline{V}_{\mathcal{H}_{1}}}l\left(Z(x^{0}(c^{*}, \cdot, 0)) + A_{1}(\cdot)y(\varepsilon) + \mathcal{R}(y(\varepsilon), \cdot, \varepsilon) \right) = 0
\end{equation}
the set of strong generalized solutions of the boundary-value problem (\ref{today:3}), (\ref{today:4}) has the following form
\begin{equation} \label{today:5} 
	y_{n}(\varepsilon) = U(n)PP_{N(D)}c(\varepsilon) + \overline{y}_{n}(\varepsilon),
\end{equation}
where 
\begin{align}
\overline{y}_{n}(\varepsilon) = \varepsilon \overline{G[Z(y_{\cdot}(\varepsilon) + x_{\cdot}^{0}(c^{*})), 0]}(n). \nonumber
\end{align}
Substituting representation (\ref{today:5}) into the generalized solvability condition (\ref{today:4}) we obtain the following operator equation
\begin{equation} \label{today:6}
	B_{0}c(\varepsilon) = -\left[ \begin{array}{ccccc} \sum_{k = -\infty}^{+\infty}\overline{H}(k + 1)\left(A_1(k)\overline{y}_{k}(\varepsilon)  + \mathcal{R}(y_{k}(\varepsilon), k, \varepsilon)\right)\\
 P_{\overline{V}_{\mathcal{H}_{1}}}l\left(A_{1}(\cdot)\overline{y}_{\cdot}(\varepsilon) + \mathcal{R}(y_{\cdot}(\varepsilon), \cdot, \varepsilon) \right)
 \end{array}
 \right]
\end{equation}  
with operator $B_{0}$ in the form 
\begin{align}
B_{0} = \left[ \begin{array}{ccccc} \sum_{k = -\infty}^{+ \infty}\overline{H}(k + 1)U(k) \\
P_{\overline{V}_{\mathcal{H}_{1}}}lA_{1}(\cdot)U(\cdot)\end{array} \right]PP_{N(D)}P_{N(V)}. \nonumber
\end{align}

Since the operator $B_{0}$ in a Hilbert space $\mathcal{H}$ always has a strong Moore-Penrose pseudoinverse $\overline{B}_{0}^{+}$ \cite{BoiPok100} the condition for the strong generalized solvability of equation (\ref{today:6}) takes the following form
\begin{align}
P_{\overline{\mathcal{H}}_{B_{0}}}\left[ \begin{array}{ccccc} \sum_{k = -\infty}^{+\infty}\overline{H}(k + 1)\left(A_1(k)\overline{y}_{k}(\varepsilon)  + \mathcal{R}(y_{k}(\varepsilon), k, \varepsilon)\right) \nonumber \\
P_{\overline{V}_{\mathcal{H}_{1}}}l\left(A_{1}(\cdot)\overline{y}_{\cdot}(\varepsilon) + \mathcal{R}(y_{\cdot}(\varepsilon), \cdot, \varepsilon) \right)\end{array} \right] = 0, \nonumber
\end{align}
where $P_{\overline{\mathcal{H}}_{B_{0}}}$ is orthoprojector onto the subspace $\overline{\mathcal{H}}_{B_{0}}$ ($\overline{\mathcal{H}}_{B_{0}} = \widetilde{\overline{\mathcal{H}}}_{\overline{D}} \ominus R(\overline{B}_{0})$, here $\widetilde{\overline{\mathcal{H}}}_{\overline{D}}$ is the completion of $\overline{\mathcal{H}}_{\overline{D}}$ according to the corresponding norm).
Since we have $\overline{H}(k + 1) = P_{\overline{\mathcal{H}}_{\overline{D}}}QU^{-1}(k + 1)$, then a sufficient condition for the strong generalized solvability of a nonlinear equation (\ref{today:3}) is the following condition
\begin{equation} \label{tod:1}
	P_{\overline{\mathcal{H}}_{B_{0}}}\left[ \begin{array}{ccccc} P_{\overline{\mathcal{H}}_{\overline{D}}}Q \\
 P_{\overline{V}_{\mathcal{H}_{1}}} 
 \end{array} \right]
 = 0.
\end{equation} 
Under condition (\ref{tod:1}) the set of strong generalized solutions of the equation (\ref{today:6}) has the following form
\begin{equation} \label{tod:2}
	c(\varepsilon) = -\overline{B}_{0}^{+}\left[ \begin{array}{ccccc}\sum_{k = -\infty}^{+\infty}\overline{H}(k + 1)\left(A_1(k)\overline{y}_{k}(\varepsilon)  + \mathcal{R}(y_{k}(\varepsilon), k, \varepsilon)\right) \\
 P_{\overline{V}_{\mathcal{H}_{1}}}l\left(A_{1}(\cdot) \overline{y}_{\cdot}(\varepsilon) + \mathcal{R}(y_{\cdot}(\varepsilon), \cdot, \varepsilon) \right)
 \end{array} \right] + P_{N(B_{0})}c_{\rho}(\varepsilon),
\end{equation}
here  $c_{\rho}(\varepsilon) \in \mathcal{H}$. Thus, the problem of the existence of a strong generalized solution $y_{n}(\varepsilon)$ of the equation (\ref{today:3}) bounded on the entire $\mathbb{Z}$ axis reduces to the strong generalized solvability of the following operator system
\begin{equation} \label{tod:3}
	\left\{ \begin{array}{ccccc}
		y_{n}(\varepsilon) = U(n)PP_{N(D)}c + \overline{y}_{n}(\varepsilon), \\
		c(\varepsilon) = -\overline{B}_{0}^{+}\left[ \begin{array}{ccccc}\sum_{k = -\infty}^{+\infty}\overline{H}(k + 1)\left(A_1(k)\overline{y}_{k}(\varepsilon)  + \mathcal{R}(y_{k}(\varepsilon), k, \varepsilon)\right) \\
 P_{\overline{V}_{\mathcal{H}_{1}}}l\left(A_{1}(\cdot) \overline{y}_{\cdot}(\varepsilon) + \mathcal{R}(y_{\cdot}(\varepsilon), \cdot, \varepsilon) \right)
 \end{array} \right] \\
		\overline{y}_{n}(\varepsilon) = \varepsilon \overline{G[Z(y_{\cdot}(\varepsilon) + x_{\cdot}^{0}(c^{*})), 0]}(n)
	\end{array}\right. + P_{N(B_{0})}c_{\rho}(\varepsilon)
\end{equation}
or in the form
\begin{equation} \label{tod:4}
	\left(\begin{array}{ccccc} y_{n}(\varepsilon) \\
		c(\varepsilon) \\
		\overline{y}_{n}(\varepsilon)\end{array}\right) = \left(\begin{array}{ccccc} 0 & U(n)PP_{N(D)} & I \\
		0 & 0 & L_{1} \\
		0 & 0 & 0 \end{array} \right) \left(\begin{array}{ccccc} y_{n}(\varepsilon) \\
		c(\varepsilon) \\
		\overline{y}_{n}(\varepsilon)\end{array}\right) + g_{n}(\varepsilon),
\end{equation}
where 
\begin{align} 
L_1* = - \overline{B}_{0}^{+}\left[ \begin{array}{ccccc}\sum_{k = -\infty}^{+\infty}\overline{H}(k + 1)A_1(k)* \nonumber \\
P_{\overline{V}_{\mathcal{H}_{1}}}lA_{1}(\cdot)*
\end{array}
\right]
\end{align}
\begin{align}
g_{n}(\varepsilon) = \left( \begin{array}{ccccc} 0 \nonumber \\  
	- \overline{B}_{0}^{+}\left[ \begin{array}{ccccc} \sum_{k = -\infty}^{+\infty}\overline{H}(k + 1)\mathcal{R}(y_{k}(\varepsilon), k, \varepsilon)  \nonumber \\
 P_{\overline{V}_{\mathcal{H}_{1}}}l\mathcal{R}(y_{\cdot}(\varepsilon), \cdot, \varepsilon)\end{array}
 \right] \nonumber \\
 \varepsilon\overline{G[Z(y_{\cdot}(\varepsilon) + x_{\cdot}^{0}(c^{*})), 0]}(n) )
\end{array} \right) + P_{N(B_{0})}c_{\rho}(\varepsilon).
\end{align}
If we denote the operator on the right hand side of the system (\ref{tod:4}) and the vector as
\begin{align}
S(n) = \left(\begin{array}{ccccc} 0 & U(n)PP_{N(D)} & I \\
	0 & 0 & L_{1} \\
	0 & 0 & 0 \end{array} \right), ~~ z_{n}(\varepsilon) = \left( \begin{array}{ccccc} y_{n}(\varepsilon) \\
	c(\varepsilon) \\
	\overline{y}_{n}(\varepsilon)
\end{array} \right), \nonumber
\end{align}
then we obtain the operator system
\begin{equation} \label{apr:1}
	(I - S(n))z_{n}(\varepsilon) = g_{n}(\varepsilon).
\end{equation} 
It is easy to see that the operator $(I - S(n))$ has a bounded inverse in the following form
\begin{equation} \label{apr:2}
	(I - S(n))^{-1} = \left( \begin{array}{cccc} 
		I & U(n)PP_{N(D)} & U(n)PP_{N(D)}L_1 + I \\
		0 & I & L_1 \\
		0 & 0 & I
	\end{array}\right).
\end{equation}
Thus, the operator system (\ref{apr:1}) can be represented as follows
\begin{equation} \label{apr:3}
	z_{n}(\varepsilon) = (I - S(n))^{-1}F_{n}(\varepsilon)z_{n}(\varepsilon),
\end{equation}
where $g_{n}(\varepsilon) = F_{n}(\varepsilon)z_{n}(\varepsilon)$. 
By choosing the parameter $\varepsilon \in [0, \varepsilon_{*}]$ and the element $c_{\rho}$, one can achieve that the operator $(I - S(n))^{-1}F_{n}(\varepsilon)$ on the right-hand side of the operator system (\ref{apr:3}) is contractive, and applying the contraction mapping principle \cite{Kolmog} , one can obtain a sufficient condition for the existence of strong generalized solutions of the nonlinear equation (\ref{today:1}).

{\bf Theorem 9.} (sufficient condition). {\it Suppose the equation (\ref{6.2.2}) admits a dichotomy on the semi-axes $\mathbb{Z}_{+}$ and $\mathbb{Z}_{-}$ with projectors $P$ and $Q$ respectively and the considered linear equation (\ref{6.2.1}), (\ref{6.2.1000})  has strong generalized bounded solutions $x_{n}^{0}(c)$ in the form (\ref{6.2.4}). Assume that 
\begin{align}
P_{\overline{\mathcal{H}}_{B_{0}}}\left[ \begin{array}{ccccc} P_{\overline{\mathcal{H}}_{\overline{D}}}Q  \\
P_{\overline{V}_{\mathcal{H}_{1}}} 
\end{array}
\right] = 0. \nonumber
\end{align}
Then for each element $c = c^{*}$ satisfying the equation for generating elements (\ref{today:2}) there are strong generalized solutions $x_{n}(\varepsilon)$ of the nonlinear boundary-value problem (\ref{today:1}), (\ref{today:2}) bounded on the entire $\mathbb{Z}$ axis, turn for $\varepsilon = 0$ into the generating solution $x_{n}^{0}(c^{*})$: $x_{n}(0) = x_{n}^{0}(c^{*})$. These solutions can be found using a convergent  iterative process for $\varepsilon \in [0, \varepsilon_{*}] \subset [0, \varepsilon_{0}]$
\begin{align}
y_{n}^{l + 1}(\varepsilon) = U(n)PP_{N(D)}c^{l + 1}(\varepsilon) + \overline{y}_{n}^{l + 1}(\varepsilon), \nonumber
\end{align}
\begin{align}
c^{l + 1}(\varepsilon) = -\overline{B}_{0}^{+}\left[ \begin{array}{ccccc} \sum_{k = -\infty}^{+\infty}\overline{H}(k + 1)\left(A_1(k)\overline{y}_{k}^{l + 1}(\varepsilon)  + \mathcal{R}(y_{k}^{l}(\varepsilon), k, \varepsilon)\right) \\
P_{\overline{V}_{\mathcal{H}_{1}}}l\left(A_{1}(\cdot)\overline{y}^{l + 1}_{\cdot}(\varepsilon) + \mathcal{R}(y^{l}_{\cdot}(\varepsilon), \cdot, \varepsilon) \right) \end{array}\right] + \mathcal{P}_{N(B_{0})}c_{\rho}(\varepsilon),
\end{align}
\begin{align}
\overline{y}_{n}(\varepsilon) = \varepsilon \overline{G[Z(y_{\cdot}(\varepsilon) + x_{\cdot}^{0}(c^{*})), 0]}(n), \nonumber
\end{align}
\begin{align}
x_{n}^{l}(\varepsilon) = y_{n}^{l}(\varepsilon) + x_{n}^{0}(c^{*}), y_{n}^{0}(\varepsilon) = 0, l = \overline{0, \infty}. \nonumber
\end{align}
}

{\bf Remark 10.} {\it It should be noted that if we choose the element $c_{\rho}(\varepsilon) \in \mathcal{H}$ in the form $c_{\rho}(\varepsilon) = \varepsilon \tilde{c}_{\rho}$ then we can always achieve that the operator $(I - S(n))^{-1}F_{n}(\varepsilon)$ is contractive.}

{\bf Remark 11.} {\it The number of strong generalized bounded solutions of the nonlinear boundary-value problem (\ref{today:1}), (\ref{today:2}) depends on the dimension of the space $N(B_{0})$. If $P_{N(B_{0})} = 0$ then we have the unique strong generalized bounded solution of the boundary-value problem (\ref{today:1}), (\ref{today:2}). If $P_{N(B_{0})} \neq 0$ then we get that new solutions $x_{n}(\varepsilon)$ appear from the point $\varepsilon = 0$ ($x_{n}(0) = x_{n}^{0}(c^{*})$). }   

{\bf Remark 12.}	{\it 	We can consider a more general boundary value problem with boundary conditions in the form 
\begin{align}
lx_{\cdot}(\varepsilon) = \alpha + \varepsilon J(x_{\cdot}(\varepsilon), \varepsilon) \nonumber
\end{align}
with a nonlinear operator-valued function $J$. We can choose among the set of bounded solutions those that are needed (for example periodic, homoclinic or heteroclinic).
} 

{\bf Theorem 10. }	{\it Suppose that operator $F(c)$ has a Frechet derivative for an element $c = c^{*}$ that satisfies the equation for generating elements (\ref{today:2}). If $F'(c^{*})$ has bounded inverse the boundary-value problem (\ref{today:1}), (\ref{today:2}) has a unique bounded solution. 
}

{\bf Proof.} It is easy to check that in this case $F'(c^{*}) = B_{0}$. It follows from this equality that the conditions of theorem are satisfied. Conditions for the invertibility of the operator $B_{0}$ connect the necessary and sufficient conditions of the existence of bounded solutions.  

{\bf Remark 13.}	{\it In the finite-dimensional case, the condition for the invertibility of the operator $B_{0}$ is the condition for the simplicity of the root $c = c^{*}$ and we can obtain a discrete analogue of the well-known Palmer's theorem  \cite[p.408]{Chuesh}, \cite{Palm1} and the Melnikov conditions of the existence of chaos (in the sense of Bernoulli). In more general case we obtain the weak conditions of chaos. }

{\bf Remark 14.} {\it In such way we obtain conditions of the existence of weak chaos and moreover the so-called reducibility conditions \cite{Hu}. We have more general case when the considered problem can be ill-posed. It means that the initial equation has the set of bounded solutions (the uniqueness of the solution is disturbed)}.

{\bf Remark 15.}{\it Existence of bounded solutions for the fractional difference equations was considered in the paper \cite{Dibl}.}

{\bf Example 1.} Consider the case when $\mathcal{H} = l_{2}$ and we have homogeneous equation $h_{n} = 0$ with operators $A_{n}$ which acts from the Hilbert space $l_{2}$ into itself by the rule ($x_{n} = (x_{n}^{1}, x_{n}^{2}, ..., x_{n}^{k}, ...) \in l_{2}$)
\begin{align}
A_{n}x_{n} = \left( 2^{-sgn(n)}x_{n}^{1}, 2^{sgn(n)}x_{n}^{2}, 2^{-sgn(n)}x_{n}^{3}, 2^{sgn(n)}x_{n}^{4}, 2^{-|sgn(n)|}x_{n}^{5}, 2^{-|sgn(n)|}x_{n}^{6}, ... \right) \nonumber
\end{align}
with boundary condition 
\begin{align}
lx_{\cdot} = x_{m} = \alpha, \alpha = (\alpha_{1}, 0, \alpha_{2}, 0, 0, ...) \in l_{2}, ~ m > 0. \nonumber
\end{align}
In this case evolution operator $U(n)$ has a diagonal form 
\begin{align}
U(n) = \left\{ \begin{array}{ccccc}
diag(2^{-(n - 1)}, 2^{n - 1}, 2^{-(n - 1)}, 2^{n - 1}, 2^{-(n - 1)}, 2^{-(n - 1)}, ...), n > 0, \\
I, ~~n = 0, \\
diag(2^{n}, 2^{-n}, 2^{n}, 2^{-n}, 2^{-n}, ...), n < 0
\end{array}\right. \nonumber
\end{align}
and corresponding equation admits an exponential dichotomy with projectors $P$ and $Q$ in the form of diagonal operators 
\begin{align}
P = diag(1, 0, 1, 0, 1, 1, ...),~~ Q = diag(0, 1, 0, 1, 1, ...). \nonumber
\end{align}
In this case operator $D: l_{2} \rightarrow l_{2}$ is selfadjoint and has a form 
\begin{align}
D = diag(0, 0, 0, 0, 0, 1, 1, 1, 1, ...) = D^{*}. \nonumber
\end{align}
Since $PQ = QP$ and $D = D^{*}$ there is an operator $D = D^{+}$ (we have that condition of corollary 5 is satisfied). Orthoprojectors $P_{N(D)} = P_{\mathcal{H}_{D}}$ has a form
\begin{align}
P_{N(D)} = I - D^{2} = P_{\mathcal{H}_{D}} = diag(1, 1, 1, 1, 0, 0, ...) \nonumber
\end{align}
and corresponding equation has a set of bounded solutions in the form
\begin{align}
x_{n}^{0}(c) = \left\{ \begin{array}{ccccc} 
(2^{-(n - 1)}c_{1}, 0, 2^{-(n - 1)}c_{2}, 0, ...), n > 0, \nonumber \\
(c_{1}, 0, c_{2}, 0, ...), n = 0, \nonumber \\
(2^{n}c_{1}, 0, 2^{n}c_{2}, 0, ...), n < 0.
\end{array} \right.
\end{align}
Using theorem 1 we obtain the operator equation 
\begin{align}
Vc = (2^{-(m - 1)}c_{1}, 0, 2^{-(m - 1)}c_{2}, 0, ...) = (\alpha_{1}, 0, \alpha_{2}, 0, ...) \nonumber
\end{align}
and unique bounded solution of the boundary-value problem (\ref{6.2.1}), (\ref{6.2.1000}) has a form
\begin{align}
x_{n}^{0} = \left\{ \begin{array}{ccccc} 
(2^{m - n}{\alpha}_{1}, 0, 2^{m - n}{\alpha}_{2}, 0, ...), n > 0, \nonumber \\
(2^{m - 1}{\alpha}_{1}, 0, 2^{m - 1}{\alpha}_{2}, 0, ...), n = 0, \nonumber \\
(2^{n + m - 1}{\alpha}_{1}, 0, 2^{n + m - 1}{\alpha}_{2}, 0, ...), n < 0.
\end{array} \right.
\end{align}

{\bf Example 2.} Let the space $\mathcal{H} = l_{2}$ and operators $A_{n}$ have the following form 
\begin{align}
A_{n}x_{n} = \left( 2^{sgn(n)}x_{n}^{1}, ..., 2^{sgn(n)}x_{n}^{k}, 2^{-sgn(n)}x_{n}^{k + 1}, 2^{-sgn(n)}x_{n}^{k + 2}, ... \right).
\end{align}
In this case evolution operator $U(n)$ has a diagonal form 
\begin{align}
U(n) = \left\{ \begin{array}{ccccc}
diag(\underbrace{2^{n - 1}, ...,  2^{n - 1}}_{k}, 2^{-(n - 1)}, 2^{-(n - 1)}, ...), n > 0, \nonumber \\
I, ~~n = 0, \nonumber \\
diag(\underbrace{2^{-n}, ..., 2^{-n}}_{k}, 2^{n}, 2^{n}, ...), n < 0
\end{array}\right.
\end{align}
and corresponding equation admits an exponential dichotomy with projectors $P$ and $Q$ in the form of diagonal operators 
\begin{align}
P = diag(\underbrace{0, ..., 0}_{k}, 0, 1, 1, ...),~~ Q = diag(\underbrace{1, ..., 1}_{k}, 0, 0, ...). \nonumber
\end{align}
In this case operator $D = 0$ and $P_{N(D)} = P_{B_{D}} = I$. Corresponding homogeneous equation has a set of bounded solutions in the form 
\begin{align}
x_{n}^{0}(c) = \left\{ \begin{array}{ccccc} 
(\underbrace{0, ..., 0}_{k}, 2^{-(n - 1)}c_{k + 1}, 2^{-(n - 1)}c_{k + 2}, ...), n > 0, \nonumber \\
(\underbrace{0, ..., 0}_{k}, c_{k + 1}, c_{k + 2}, ...), n = 0, \nonumber \\
(\underbrace{0, ..., 0}_{k}, 2^{n}c_{k + 1}, 2^{n}c_{k + 2}, ...), n < 0.
\end{array} \right. \nonumber
\end{align}
Nonhomogeneous equation (\ref{6.2.1}) ($h_{n} \neq 0$) has bounded solutions if and only if the following $k$ conditions are satisfied
\begin{align}
\sum_{l = -\infty}^{- 1}2^{l + 1}h_{l}^{p} + h_{0}^{p} + \sum_{l = 1}^{+\infty}2^{-l}h_{l}^{p} = 0, p = \overline{1, k} \nonumber
\end{align}
and in this case generalized Green's operator has the following form $(G[h])(n) = $
\begin{equation} \label{dubov:1}
= \left\{ \begin{array}{cccccc}
\left( -\sum_{l = n}^{+\infty}2^{n - l - 1}h_{l}^{1}, ..., -\sum_{l = n}^{+\infty}2^{n - l - 1}h_{l}^{k}, \sum_{l = 0}^{n - 1}2^{- n + l + 1}h_{l}^{k + 1}, ... \right), n \geq 0, \\
\left(\sum_{l = -\infty}^{n - 1}2^{-n + l + 1}h_{l}^{1}, ..., \sum_{l = \infty}^{n - 1}2^{- n + l + 1}h_{l}^{k}, \sum_{l = n}^{-1}2^{n - l - 1}h_{l}^{k + 1}, ... \right), n \leq 0. 
\end{array} \right.
\end{equation}

The set of bounded solutions of the equation (\ref{6.2.1}) has the form 
\begin{align}
x_{n}^{0}(c) = \overline{x}_{n}^{0}(c) + (G[h])(n). \nonumber
\end{align}

{\bf Example 3.} Consider homogeneous equation ($h_{n} = 0$) from the example 2 with boundary condition 
\begin{equation} \label{dub:2}
lx_{\cdot} = A(x_{q} - x_{p}) = (\underbrace{0, ..., 0}_{k}, x_{q}^{k + 2} - x_{p}^{k + 2}, ..., \frac{x_{q}^{k + m} - x_{p}^{k + m}}{m}, ...) = \alpha,
\end{equation}
where $\alpha = (\underbrace{0, ...., 0}_{k + 1}, \alpha_{k + 2}, \alpha_{k + 3}, ...) \in l_{2}, 0 < q < p$. Substituting solutions $\overline{x}_{n}^{0}(c)$ in boundary condition (\ref{dub:2}) we obtain the following operator equation
\begin{align}
Vc = (2^{- q + 1} - 2^{ - p + 1})\left(\underbrace{0, ..., 0}_{k + 1}, c_{k + 2}, \frac{c_{k + 3}}{2}, ..., \frac{c_{k + m + 1}}{m}, ...  \right) =
\nonumber
\end{align}
\begin{align}
= (\underbrace{0, ..., 0}_{k + 1}, \alpha_{k + 2}, \alpha_{k + 3}, ...).
\nonumber
\end{align}

It is easy to check that the operator $V$ has a nonclosed set of values. Really, the following sequence of elements $c_{n} = (\underbrace{0, ..., 0}_{k}, \underbrace{1, ..., 1}_{n}, 0, ...) \in l_{2}$ belongs to the set of values of $V$, but the limit element $c = (\underbrace{0, ..., 0}_{k}, 1, ..., 1, ...)$ does not belong to the set of values of $V$. Kernel $N(V)$ of the operator $V$ consists from the elements $\{(0, ..., 0, c, 0, ...), c \in \mathbb{R} \} $. After factorization we have the operator 
\begin{align}
\mathcal{V} : X \rightarrow R(V), \mathcal{V}c = \mathcal{V}(0, ..., 0, c_{k + 2}, c_{k + 3}, ...) = (0, ..., 0, c_{k + 2}, \frac{c_{k + 3}}{2}, ..., \frac{c_{k + m + 1}}{m}, ...)
\nonumber
\end{align}
which is linear, injective and continuous. Completing the space $X$ according to the norm 
\begin{align}
||c||_{\overline{X}}^{2} = ||\mathcal{V}c||^{2}_{R(V)} = \sum_{i = 1}^{+\infty}\frac{c_{k + i + 1}^{2}}{i^{2}}
\nonumber
\end{align}
we obtain a new space $\overline{X}$ and $l_{2} \subset \overline{\mathcal{H}} = N(V) \oplus \overline{X}$ (see \cite{BoiPok100}). It is easy to check that in this case operator $V$ has a strong Moore-Penrose pseudo-inverse in the following form 
\begin{align}
\overline{V}^{+} \alpha = \frac{1}{2^{- q + 1} - 2^{- p + 1}}(\underbrace{0, ..., 0}_{k + 1}, \alpha_{k + 2}, 2\alpha_{k + 3}, ..., m\alpha_{k + m + 1}, ...).
\nonumber
\end{align}

Thus, we obtain that the boundary-value problem (\ref{6.2.1}), (\ref{6.2.1000}) has a set of strong generalised bounded solutions in the form
\begin{align}
x_{n}^{0}(c) = \left\{ \begin{array}{ccccc}(\underbrace{0, ..., 0}_{k}, 2^{-(n - 1)}c, \frac{2^{-n}}{2^{-q} - 2^{- p}}\alpha_{k + 2}, ..., \frac{2^{-n}m}{2^{-q} - 2^{-p}}\alpha_{k + m + 1}, ...),  n > 0, \\
(\underbrace{0, ..., 0}_{k}, c, \frac{1}{2^{-q + 1} - 2^{- p + 1}}\alpha_{k + 2}, ..., \frac{m}{2^{-q + 1} - 2^{-p + 1}}\alpha_{k + m + 1}, ...),  n = 0, \\
(\underbrace{0, ..., 0}_{k}, 2^{n}c, \frac{2^{n}}{2^{-q + 1} - 2^{- p + 1}}\alpha_{k + 2}, ..., \frac{2^{n}m}{2^{-q + 1} - 2^{-p + 1}}\alpha_{k + m + 1}, ...),  n < 0, \\
\end{array}
\right.
\nonumber
\end{align}
for any $c \in \mathbb{R}$.

It should be noted that in this work author use well-known results which was described in the monograph \cite{BoiSam1}. With using technique of generalized solutions which was developed in the following monograph \cite{Lyash} author can developed a new method for investigating of weak chaos in discrete systems. Such theory can be applied to investigation of the problems to which the following works are devoted \cite{Nonl1}, \cite{Nonl2}, \cite{Nonl3}, \cite{Nonl4}, \cite{Nonl5} (see also the following works \cite{Nonl6}, \cite{Hopf}).

{\bf Conclusions.} Proposed in the article statements gives us possibility to investigate the question of the existence of bounded solutions of the linear and weakly nonlinear nonhomogeneous equation in the Banach and Hilbert spaces in general. Developed in the paper method
allows us to investigate boundary value problems on the whole axis with conditions at infinity. As application we can consider the countable system of difference equations. Obtained in the work the necessary condition of solvability is an analogue of Fredholm's alternative.

\end{document}